\documentclass[11pt]{article}
\usepackage{amssymb, amsthm, amsmath}

\newtheorem{thm}{Theorem}[section]
\newtheorem{lem}[thm]{Lemma}
\newtheorem{prop}[thm]{Proposition}

\newcommand{\N}{\mathbb{N}}
\newcommand{\Z}{\mathbb{Z}}
\newcommand{\R}{\mathbb{R}}
\newcommand{\T}{\mathbb{T}}

\newcommand{\Punc}{\operatorname{Punc}}

\title{Affability of equivalence relations \\
arising from two-dimensional substitution tilings}
\author{MATUI Hiroki}
\date{}

\begin{document}
\maketitle

\begin{abstract}
We will show that an equivalence relation on a Cantor set 
arising from a two-dimensional substitution tiling by polygons 
is affable in the sense of Giordano, Putnam and Skau. 
\end{abstract}

\section{Introduction}

We study an equivalence relation arising from 
a substitution tiling system in $\R^2$. 
The dynamics of substitution tiling systems have been studied 
by many authors (see \cite{KP} and the references given there). 
A tiling in $\R^2$ gives rise to an action of $\R^2$ 
on a topological space $\Omega$. 
Each element of this space can be thought of as a tiling of $\R^2$ 
and the $\R^2$-action is given by translation. 
Under suitable hypotheses, 
this topological space $\Omega$ becomes compact 
and the action of $\R^2$ becomes minimal. 
Besides, we can find a transversal to this $\R^2$-action 
and it is known to be a Cantor set. 
Our object of study is an equivalence relation on this Cantor set 
which is induced from the $\R^2$-action. 
We remark that $\Omega$ is known to be homeomorphic to 
the suspension of a free minimal $\Z^2$-action 
on a Cantor set $X$ (see \cite{SW}). 
This means that 
the equivalence relation arising from a tiling space is isomorphic to 
an induced equivalence relation on a certain clopen subset $U\subset X$. 

The study of topological orbit equivalence 
between equivalence relations on a Cantor set was initiated 
by Giordano, Putnam and Skau in \cite{GPS1}. 
They showed that equivalence relations arising from $\Z$-actions 
on Cantor sets are (topologically) orbit equivalent to 
AF equivalence relations. 
An equivalence relation is said to be AF 
when it can be written as an inductive limit of 
compact open subequivalence relations. 
An equivalence relation is said to be affable 
when it is orbit equivalent to an AF equivalence relation (\cite{GPS2}). 
In a recent paper \cite{GPS3}, they also showed that 
equivalence relations arising from $\Z^2$-actions on Cantor sets 
are affable under a hypothesis involving the existence of 
sufficiently many small positive cocycles. 

In the case of tilings, however, 
it is not known whether the equivalence relation admits 
sufficiently many such cocycles. 
It relates to a more complete understanding of the cohomology 
of tiling spaces. 
In this paper we take a more direct tack. 
We will prove the affability by a similar method to \cite{M}, 
in which any extension of a product of 
two Cantor minimal $\Z$-systems was shown to be affable. 
As discussed in \cite{P}, 
the equivalence relation arising from a substitution tiling has 
a canonical AF subequivalence relation. 
But, this AF equivalence relation is too large to apply 
the absorption theorem in \cite{GPS2}. 
We will construct a smaller AF subequivalence relation carefully 
so that the absorption theorem can be applied.

\section{Preliminaries}

Let ${\cal V}$ be a finite collection of polygons in $\R^2$. 
Each element of ${\cal V}$ is called a prototile. 
Every prototile is homeomorphic to the closed unit ball 
and has finitely many edges and vertices. 
We call a translate of one of the prototile a tile. 
Let $\omega$ be a substitution rule and 
let $\lambda>1$ be its inflation constant. 
Thus, for each prototile $p\in{\cal V}$, 
$\omega(p)$ is a finite collection of tiles 
with pairwise disjoint interiors and 
their union is $\lambda p=\{\lambda v:v\in p\}$. 
Such a pair $({\cal V},\omega)$ is called 
a (two-dimensional) substitution tiling system. 
Several examples, including the Penrose tiles, were given in \cite{AP}. 
We also allow the possibility that our prototiles carry labels. 
Even if two tiles coincide as geometric objects, 
they may be distinguished. 
A tiling is a collection of tiles such that their union is $\R^2$ and 
their interiors are pairwise disjoint. 

As in \cite{AP}, \cite{KP} and \cite{P}, 
we assume that our substitution tiling system $({\cal V},\omega)$ is 
primitive, aperiodic and satisfies the finite pattern condition. 
We follow the notation of \cite{P}. 
Thanks to the finite pattern condition, 
we can define a compact metric space $\Omega$ consisting of tilings. 
For each tiling $T\in\Omega$, 
we let $\omega(T)=\{\omega(p):p\in T\}$. 
By the aperiodicity condition, 
this map $\omega:\Omega\to\Omega$ becomes a homeomorphism. 
We will also assume that our substitution system forces its border. 
As discussed in \cite{KP}, this loses no generality, 
provided we allow labelled tiles. 

Without loss of generality, we may assume that 
every $p\in{\cal V}$ contains the origin in its interior. 
We define the puncture of any prototile $p\in{\cal V}$ 
to be the origin. 
For any $v\in\R^2$, 
we define the puncture of the tile $p+v$ to be $v$. 
As in \cite{P}, 
$\Omega_{punc}$ denotes the set of all tilings $T$ in $\Omega$ 
such that the origin is a puncture of some tile $p\in T$. 
This set is compact, totally disconnected and has no isolated points, 
that is, $\Omega_{punc}$ is a Cantor set. 

Define an equivalence relation $R_{punc}$ on $\Omega_{punc}$ by 
\[ R_{punc}=\{(T,T+v):T,T+v\in\Omega_{punc}, \ v\in\R^2\}. \]
As explained in \cite{P}, 
this equivalence relation is given a topology 
so that it becomes a locally compact, Hausdorff, $\sigma$-compact, 
r-discrete, principal groupoid. 
In other words, $R_{punc}$ is an \'etale equivalence relation 
\cite[Definition 2.1]{GPS2}. 
Since the action of $\R^2$ on $\Omega$ by translation is minimal, 
the equivalence relation $R_{punc}$ becomes minimal. 
In \cite{P}, a subequivalence relation $R_{AF}\subset R_{punc}$ 
was introduced. 
The equivalence relation $R_{AF}$ is a minimal AF equivalence relation 
with the relative topology from $R_{punc}$. 
Our main theorem asserts that 
$R_{punc}$ is orbit equivalent to $R_{AF}$. 
\bigskip

For each $p\in{\cal V}$, let $\Punc(p)$ denote 
the set of all punctures in the tiles of $\omega(p)$. 
Let ${\cal E}$ be the disjoint union of 
all the $\Punc(p)$'s for $p\in{\cal V}$. 
We define $s:{\cal E}\to{\cal V}$ and $r:{\cal E}\to{\cal V}$ as follows: 
When $e\in{\cal E}$ is contained in $\Punc(q)$, we put $r(e)=q$. 
When $p\in{\cal V}$ satisfies $p+e\in\omega(r(e))$, we put $s(e)=p$. 
We regard $({\cal V},{\cal E},r,s)$ as a finite directed graph. 
Since the substitution is primitive, this graph is also primitive. 

Put 
\[ X=\{(e_n)_{n=1}^\infty\in{\cal E}^\N:
r(e_n)=s(e_{n+1})\text{ for all }n\in\N\}. \]
It is easy to see that $X$ is a Cantor set 
with the induced topology from ${\cal E}^\N$. 
Moreover, there exists a natural homeomorphism 
$\pi:\Omega_{punc}\to X$ defined as follows: 
For $T\in\Omega_{punc}$ and $n\in\N\cup\{0\}$, 
let $p_n+v_n$ be the tile in $\omega^{-n}(T)$ 
whose interior contains the origin, where $p_n\in{\cal V}$ and $v_n\in\R^2$. 
Note that $v_0$ equals zero because $T$ is in $\Omega_{punc}$. 
For every $n\in\N$, it is not hard to see that 
$e_n=v_{n-1}-\lambda v_n$ belongs to $\Punc(p_n)$. 
By defining $\pi(T)=(e_n)_n$, 
we obtain a homeomorphism $\pi:\Omega_{punc}\to X$. 

We denote the copies of $R_{punc}$ and $R_{AF}$ 
via the map $\pi:\Omega_{punc}\to X$ 
by ${\cal X}_\omega$ and ${\cal X}$, respectively. 
Thus, 
\[ {\cal X}_\omega=\{(x,y)\in X\times X:
(\pi^{-1}(x),\pi^{-1}(y))\in R_{punc}\} \]
and 
\[ {\cal X}=\{(x,y)\in X\times X:
(\pi^{-1}(x),\pi^{-1}(y))\in R_{AF}\}. \]
We topologize ${\cal X}_\omega$ and ${\cal X}$ 
by transferring the topology of $R_{punc}$ and $R_{AF}$ 
via the map $\pi\times\pi$. 
The aim of this paper it to show that ${\cal X}_\omega$ is affable. 
For every $n\in\N$, we define a subequivalence relation 
${\cal X}_n\subset {\cal X}$ by 
\[ {\cal X}_n=\{((e_k)_k,(f_k)_k)\in X\times X:
e_m=f_m\text{ for all }m>n\}. \]
Clearly ${\cal X}_n$ is a compact open subequivalence relation 
of ${\cal X}$ and ${\cal X}=\bigcup_{n\in\N}{\cal X}_n$. 
\bigskip

For a prototile $p\in{\cal V}$ 
we denote the set of all edges of $p$ by $E(p)$. 
Let $v_1,v_2,\dots,v_n\in\R^2$ be the vertices of $p$, 
which are ordered counterclockwise. 
When $a_i\in E(p)$ is the edge between $v_i$ and $v_{i+1}$, 
we define 
\[ \theta_p(a_i)
=\frac{v_{i+1}-v_i}{\lvert v_{i+1}-v_i\rvert}\in\T, \]
where the indices are understood modulo $n$ and 
$\T$ denotes the set of unit vectors in $\R^2$. 

We need the following hypotheses about edges of prototiles. 
Let $p,q\in{\cal V}$ be prototiles 
and suppose $p+v\in\omega(q)$ for $v\in\R^2$. 
\begin{enumerate}
\item[(P1)] If $a\in E(p)$ is an edge between vertices $u_1$ and $u_2$, 
then $\lambda^{-1}a$ is contained in the boundary of $q-\lambda^{-1}v$, 
or $\lambda^{-1}(a\setminus\{u_1,u_2\})$ is contained 
in the interior of $q-\lambda^{-1}v$. 
\item[(P2)] Let $a_1,a_2\in E(p)$ with $\theta_p(a_1)=\theta_p(a_2)$. 
If there exist edges $b_1,b_2\in E(q)$ such that 
$\lambda^{-1}(a_i+v)\subset b_i$ for $i=1,2$, 
then we have $a_1=a_2$ and $b_1=b_2$. 
\item[(P3)] Let $a_1,a_2\in E(p)$ with $\theta_p(a_1)\neq\theta_p(a_2)$. 
If there exist edges $b_1,b_2\in E(q)$ such that 
$\lambda^{-1}(a_i+v)\subset b_i$ for $i=1,2$, 
then $a_1$ is adjacent to $a_2$ and $b_1$ is adjacent to $b_2$. 
\end{enumerate}
We further assume the following. 
\begin{enumerate}
\item[(P4)] For any prototiles $p_1,p_2\in{\cal V}$ and 
their edges $a_i\in E(p_i)$, the length of $\lambda a_1$ is greater 
than the length of $2a_2$. 
\end{enumerate}
We remark that (P4) can be achieved 
by replacing $\omega$ and $\lambda$ with $\omega^N$ and $\lambda^N$ 
for sufficiently large $N\in\N$. 

Let $e\in{\cal E}$. 
By definition, $e\in\Punc(r(e))$ and $s(e)+e\in\omega(r(e))$. 
We define 
\[ A(e)=\{(a,b):a\in E(s(e)),b\in E(r(e))\text{ and }
\lambda^{-1}(a+e)\subset b\}. \]
Notice that $\theta_q(b)$ equals $\theta_p(a)$ for every $(a,b)\in A(e)$. 
The conditions (P2) and (P3) above can be interpreted as follows: 
\begin{enumerate}
\item[(P2)] If $(a_i,b_i)\in A(e)$ for $i=1,2$ and 
$\theta_{s(e)}(a_1)=\theta_{s(e)}(a_2)$, then $a_1=a_2$ and $b_1=b_2$. 
\item[(P3)] If $(a_i,b_i)\in A(e)$ for $i=1,2$ and 
$\theta_{s(e)}(a_1)\neq\theta_{s(e)}(a_2)$, then 
$a_1$ is adjacent to $a_2$ and $b_1$ is adjacent to $b_2$. 
\end{enumerate}
By (P2) and (P3), we can see that 
the cardinality of $A(e)$ is not greater than two. 
Moreover, (P3) and (P4) imply the following. 
\begin{enumerate}
\item[(P5)] For any $q\in{\cal V}$ and $b\in E(q)$, 
there exists $e\in{\cal E}$ such that 
$r(e)=q$ and $A(e)=\{(a,b)\}$ for some $a\in E(s(e))$. 
\end{enumerate}
\bigskip

Let $x=(x_n)_n\in X$. 
We say that $a=(a_n)_n$ is a border of $x$ 
when $(a_n,a_{n+1})\in A(x_n)$ for every $n\in\N$. 
Note that if $a$ is a border of $x$, 
then $\theta_{s(x_n)}(a_n)=\theta_{s(x_m)}(a_m)$ for all $n,m\in\N$. 
We write this value by $\theta(a)$. 

For $t\in\T$, we define 
\[ B_t=\{x\in X:x\text{ has a border }a=(a_n)_n
\text{ with }\theta(a)=t\}. \]
The following proposition will play a crucial role in Section 5. 

\begin{prop}\label{etale}
For $t\in\T$, $B_t$ is a closed ${\cal X}$-\'etale subset. 
\end{prop}
\begin{proof}
Let $U_m$ be the set of $x=(x_n)_n\in X$ for which 
there exist $a_i\in E(s(x_i))$ for $i=1,2,\dots,m$ 
such that $\theta_{s(x_1)}(a_1)=t$ and 
$(a_i,a_{i+1})\in A(x_i)$ for every $i=1,2,\dots,m-1$. 
It can be easily seen that $U_m$ is clopen. 
Since $B_t$ is the intersection of all the $U_m$'s, 
it is a closed subset. 

In order to show that $B_t$ is ${\cal X}$-\'etale, 
take $x=(x_n)_n\in B_t$ and $y=(y_n)_n\in B_t$ 
with $(x,y)\in {\cal X}_m$. 
Assume that $a=(a_n)_n$ is a border of $x$ and 
$b=(b_n)_n$ is a border of $y$ 
such that $\theta(a)=\theta(b)=t$. 
Define 
\[ U=\{(z,w)\in {\cal X}_m:
x_i=z_i\text{ and }y_i=w_i\text{ for all }i=1,2,\dots,m+1\}. \]
Then, $U$ is a clopen neighborhood of $(x,y)$ in ${\cal X}$. 
Suppose that $(z,w)\in U$ and $z\in B_t$. 
We would like to show that $w$ also belongs to $B_t$. 
By the definition of $B_t$, 
$z$ has a border $c=(c_n)_n$ with $\theta(c)=t$. 
Since $(b_{m+1},b_{m+2})\in A(y_{m+1})$, 
$(c_{m+1},c_{m+2})\in A(z_{m+1})$, 
$\theta(b)=\theta(c)=t$ and $y_{m+1}=z_{m+1}$, 
thanks to the condition (P2), 
we get $b_{m+1}=c_{m+1}$ and $b_{m+2}=c_{m+2}$. 
Put 
\[ d_n=\begin{cases}
b_n & n\leq m+1 \\
c_n & n>m+1. \end{cases} \]
One can see that $d=(d_n)_n$ is a border of $w$ 
with $\theta(d)=t$, which means $w\in B_t$. 
It follows that $B_t$ is ${\cal X}$-\'etale. 
\end{proof}

\section{Boundaries of tilings}

For a tiling $T\in\Omega$, 
we let $\partial(T)$ be the union of boundaries of 
all tiles belonging to $T$. 
We define 
\[ \partial_\infty(T)
=\bigcap_{n\in\N}\lambda^n\partial(\omega^{-n}(T)). \]
It is clear that for each connected component $O$ of 
$\R^2\setminus\partial_\infty(T)$ the set 
$\{T-v\in\Omega_{punc}:v\in O\}$ is a $R_{AF}$-orbit. 
We have the following three possibilities: 
When $\partial_\infty(T)$ is an empty set, 
we say that $T$ is of type I. 
In this case the $R_{punc}$-orbit of $T$ is equal to 
the $R_{AF}$-orbit of it. 
When $\partial_\infty(T)$ is a line, 
we say that $T$ is of type II. 
In this case the $R_{punc}$-orbit of $T$ 
splits into two $R_{AF}$-orbits. 
In the other case, 
we say that $T$ is of type III. 
If $T$ is of type III, then there exists $v\in\R^2$ such that 
$\partial_\infty(T)$ equals 
a union of finitely many (at least two) half lines 
with the same end-point $v\in\R^2$. 
We call $v$ the center of $\partial_\infty(T)$. 
Note that if $\partial_\infty(T)$ consists of $n$ half lines 
then the $R_{punc}$-orbit of $T$ splits 
into $n$ $R_{AF}$-orbits. 
Evidently the type of $T$ depends 
only on the $R_{punc}$-orbit of $T$. 
We also remark that a connected component of 
$\R^2\setminus\partial_\infty(T)$ can be an open half-plane, 
even if $T$ is of type III. 

For distinct $s,t\in\T$, 
we would like to define a closed subset $B_{s,t}$ of 
$B_s\cap B_t$ as follows: 
Let $x=(x_n)_n\in B_s\cap B_t$. 
By definition, $x$ has borders $a=(a_n)_n$ and $b=(b_n)_n$ 
with $\theta(a)=s$ and $\theta(b)=t$. 
Let $v_1,v_2,\dots,v_n$ be the vertices of $s(x_1)$, 
which are ordered counterclockwise. 
Assume that $a_1$ is the edge between $v_{i-1}$ and $v_i$. 
From the condition (P3), $b_1$ is adjacent to $a_1$. 
We define $B_{s,t}$ to be the set of 
all $x\in B_s\cap B_t$ for which 
the edge $b_1$ agrees with the edge between $v_i$ and $v_{i+1}$. 
Clearly $B_{s,t}\cup B_{t,s}=B_s\cap B_t$. 

For $u=(u_1,u_2),v=(v_1,v_2)\in\R^2$, 
we denote the real number $u_1v_2-u_2v_1$ by $\det(u,v)$. 
We can see the following lemma 
by an easy geometric observation and the condition (P3). 

\begin{lem}\label{cornersurvive}
Let $x=(x_n)_n,y=(y_n)_n\in X$. 
Let $s,t\in\T$ be distinct elements 
satisfying $\det(s,t)>0$. 
Suppose that both $x$ and $y$ belong to $B_{s,t}$. 
If $r(x_n)=r(y_n)$, then we get $x_n=y_n$. 
\end{lem}

By the lemma above, we obtain the following. 

\begin{lem}\label{corner}
Let $s,t\in\T$ be distinct elements 
satisfying $\det(s,t)>0$. 
\begin{enumerate}
\item The cardinality of $B_{s,t}$ is not greater 
than the cardinality of ${\cal V}$. 
\item Every $x\in B_{s,t}$ is a periodic sequence. 
\end{enumerate}
\end{lem}
\begin{proof}
(1) Let $N$ be the cardinality of ${\cal V}$. 
Suppose that $x^{(0)},x^{(1)},\dots,x^{(N)}$ belong to $B_{s,t}$. 
There exist distinct $i$ and $j$ such that 
$\{n\in\N:r(x^{(i)}_n)=r(x^{(j)}_n)\}$ is infinite. 
It follows from the lemma above that $x^{(i)}=x^{(j)}$, 
which means that 
the cardinality of $B_{s,t}$ is not greater than $N$. 

(2) Let $x\in B_{s,t}$. 
It is obvious that $x^{(m)}=(x_{n+m})_n$ also belongs to 
$B_{s,t}$ for every $m\in\N$. 
From (1), $x$ must be periodic. 
\end{proof}

\begin{lem}\label{finiteIII}
The set $L=\{T\in\Omega_{punc}:T\text{ is of type III }\}$ 
consists of finitely many $R_{punc}$-orbits. 
\end{lem}
\begin{proof}
Suppose that $T\in\Omega_{punc}$ is of type III. 
There exists at least one connected component $O$ 
of $\R^2\setminus\partial_\infty(T)$ which is congruent to a cone, 
that is, $O$ is a translate of 
\[ \{u\in\R^2:\det(s,u)>0, \ \det(t,u)>0\} \]
for some distinct $s,t\in\T$. 
By replacing $(s,t)$ with $(t,s)$ if necessary, 
we may assume $\det(s,t)>0$. 
Choose a tile $p\in T$ whose puncture $v$ is contained in $O$. 
Put $x=(x_n)_n=\pi(T-v)$. 
It is not hard to see that $(x_{n+m})_n$ belongs to $B_{s,t}$ 
for sufficiently large $m\in\N$. 
Since $\{t\in\T:B_t\neq\emptyset\}$ is a finite set, 
it suffices to show the following: 
For any distinct $s,t\in\T$ satisfying $\det(s,t)>0$, 
\[ C=\{x=(x_n)_n\in X:\text{there exists }m\in\N
\text{ such that }(x_{n+m})_m\in B_{s,t}\} \]
consists of finitely many ${\cal X}$-orbits. 
Assume that $x^{(0)},x^{(1)},\dots,x^{(N)}$ belong to $C$, 
where $N$ is the cardinality of ${\cal V}$. 
There exists $m\in\N$ such that 
$y^{(i)}=(x^{(i)}_{n+m})_n$ is contained in $B_{s,t}$ 
for all $i=0,1,\dots,N$. 
By Lemma \ref{corner} (1), there must exist distinct $i$ and $j$ 
such that $y^{(i)}=y^{(j)}$, 
which implies that $(x^{(i)},x^{(j)})\in{\cal X}_m\subset{\cal X}$. 
Therefore $C$ consists of at most $N$ ${\cal X}$-orbits. 
\end{proof}

\begin{lem}
Suppose that $T\in\Omega$ is of type III and 
the center of $\partial_\infty(T)$ is the origin. 
Let $p\in T$ be a tile which meets the origin and 
let $w\in\R^2$ be the puncture of $p$. 
If $\omega^{-k}(T)=\omega(T)$, then 
the sequence $\pi(T-w)$ has period $k$. 
\end{lem}
\begin{proof}
For $n\in\N\cup\{0\}$, let $q_n+u_n\in\omega^{-n}(T-w)$ be the tile 
which contains the origin in its interior, 
where $q_n\in{\cal V}$ and $u_n\in\R^2$. 
It suffices to show that 
the sequence $(u_{n-1}-\lambda u_n)_n$ has period $k$. 
As $p$ meets the origin and the origin is in $\partial_\infty(T)$, 
we can see that $-\lambda^{-n}w$ is 
in the boundary of the tile $q_n+u_n\in\omega^{-n}(T-w)$. 
Hence the boundary of 
\[ q_n+u_n+\lambda^{-n}w\in\omega^{-n}(T) \]
contains the origin. Put $p_n=q_n+u_n+\lambda^{-n}w$. 
By $\omega^{-k}(T)=\omega(T)$, we have $p_n\in\omega^{-n-k}(T)$, 
and so the interior of $p_n$ does not intersect $\lambda^k p_n$ or 
$p_n$ is contained in $\lambda^k p_n$. 
We obtain $p_n\subset\lambda^k p_n$, 
because the boundary of $p_n$ contains the origin and 
$p_n$ is a polygon. 
It follows that 
\[ \lambda^{-k}(q_n+u_n)\subset
q_n+u_n+\lambda^{-n}w-\lambda^{-n-k}w. \]
From this and 
\begin{align*}
\omega^{-n-k}(T-w)
&= \omega^{-n}(T)-\lambda^{-n-k}w \\
&= \omega^{-n}(T-w)+\lambda^{-n}w-\lambda^{-n-k}w \\
&\ni q_n+u_n+\lambda^{-n}w-\lambda^{-n-k}w, 
\end{align*}
we get $q_{n+k}=q_n$ and 
$u_{n+k}=u_n+\lambda^{-n}w-\lambda^{-n-k}w$. 
Then we have 
\begin{align*}
& u_{n+k-1}-\lambda u_{n+k} \\
&= u_{n-1}+\lambda^{-n+1}w-\lambda^{-n-k+1}w
-\lambda(u_n+\lambda^{-n}w-\lambda^{-n-k}w) \\
&= u_{n-1}-\lambda u_n, 
\end{align*}
thereby completing the proof. 
\end{proof}

\begin{lem}\label{periodicIII}
When $T\in\Omega_{punc}$ is of type III, 
there exists a unique periodic sequence $x\in [\pi(T)]_{\cal X}$. 
\end{lem}
\begin{proof}
Since the uniqueness is clear, it suffices to show the existence. 

As in the proof of Lemma \ref{finiteIII}, 
we can find a connected component $O$ of 
$\R^2\setminus\partial_\infty(T)$ which is congruent to a cone, 
that is, $O$ is a translate of 
$\{u\in\R^2:\det(s,u)>0, \ \det(t,u)>0\}$ 
for some distinct $s,t\in\T$. 
By replacing $(s,t)$ with $(t,s)$ if necessary, 
we may assume $\det(s,t)>0$. 
Choose $u\in O$ so that $S=T-u\in\Omega_{punc}$. 
At first we would like to show the assertion for $S$. 
Let $y=(y_n)_n=\pi(S)$. 
There exists $m\in\N$ such that $(y_{n+m})_n$ belongs to $B_{s,t}$. 
By virtue of Lemma \ref{corner} (2), 
we see that $(y_{n+m})_n$ is a periodic sequence. 
Let $k\in\N$ be its period. 
Take $l\in\N$ so that $kl>m$. 
Define $z=(z_n)_n\in X$ by $z_n=y_{n+kl}$. 
It is easy to see that $z$ is a periodic sequence and 
$(z,y)\in{\cal X}$. 
Thus the assertion has been shown for $S$. 

Let $v\in\R^2$ be the center of $\partial_\infty(T)$. 
The sequence $z=(z_n)_n$ has period $k$, and so we have 
\[ \{p\in T-v:p\subset(\overline{O}-v)\neq\emptyset\}
=\{p\in\omega^{-k}(T-v):p\subset(\overline{O}-v)\neq\emptyset\}, \]
where $\overline{O}$ means the closure of $O$. 
Since the substitution system forces its border, 
we can conclude that $T-v$ equals $\omega^{-k}(T-v)$. 

Let $U$ be the connected component of 
$\R^2\setminus\partial_\infty(T)$ which contains the origin. 
We notice that $U$ corresponds to the ${\cal X}$-orbit 
of $\pi(T)$. 
We can find a tile $p\in T$ which meets $v$ and 
whose interior is contained in $U$. 
Then we have $\pi(T-w)\in[\pi(T)]_{\cal X}$, 
where $w$ is the puncture of $p$. 
The conclusion follows from the lemma above. 
\end{proof}

By Lemma \ref{finiteIII}, 
\[ L=\{T\in\Omega_{punc}:T\text{ is of type III}\} \]
consists of finitely many $R_{AF}$-orbits. 
It follows from Lemma \ref{periodicIII} that 
\[ M=\{x\in X:x\text{ is periodic and }\pi^{-1}(x)
\text{ is of type III}\} \]
is a finite set and 
there exists a natural bijective correspondence 
between $M$ and $\pi(L)/{\cal X}$. 
By replacing $\omega$ with $\omega^N$ for some appropriate $N\in\N$, 
we may assume that every $x\in M$ has period one. 
Define 
\[ {\cal C}=\{e\in{\cal E}:x=(e,e,\dots)\in\pi(L)\}. \]
We remark that there exists a natural bijective correspondence 
between ${\cal C}$ and $\pi(L)/{\cal X}$. 

By the conditions (P2) and (P3), we also have the following. 

\begin{lem}\label{2possibilities}
For $e\in{\cal C}$, we have $A(e)=\{(a,a)\}$ or 
$A(e)=\{(a,a),(b,b)\}$ for some $a,b\in E(r(e))$. 
\end{lem}

In the next section we will need the following condition. 
\begin{enumerate}
\item[(P6)] For any $p\in{\cal V}$ and $a\in E(p)$, 
there exists $e\in{\cal E}\setminus{\cal C}$ 
such that $A(e)=\{(b,a)\}$ for some $b\in E(s(e))$. 
\end{enumerate}
Replacing $\omega$ with $\omega^2$ if necessary, 
we can easily achieve this condition.

\section{Construction of a subequivalence relation}

In this section we would like to construct a subequivalence relation 
${\cal X}'$ of ${\cal X}$. 

At first we define a closed subset $\delta_e$ of $r(e)\in{\cal V}$ 
for each $e\in{\cal C}$ as follows. 
Let $v_1,v_2,\dots,v_n\in\R^2$ be the vertices of $r(e)$, 
which are ordered counterclockwise. 
By definition, $\lambda^{-1}(s(e)+e)$ is a subset of $r(e)$, 
and $s(e)$ is equal to $r(e)$. 
By Lemma \ref{2possibilities}, 
we have $A(e)=\{(a,a)\}$ or $A(e)=\{(a,a),(b,b)\}$ 
for some $a,b\in E(r(e))$. 

Let us consider the case of $A(e)=\{(a,a)\}$. 
Suppose that $a$ is the edge between $v_i$ and $v_{i+1}$, 
where the indices are understood modulo $n$. 
Since $\lambda^{-1}(a+e)\subset a$, 
the two points $\lambda^{-1}(v_i+e)$ and $\lambda^{-1}(v_{i+1}+e)$ 
lie on the edge $a$. 
In this case we let $\delta_e$ be the closed interval 
between $v_i$ and $\lambda^{-1}(v_{i+1}+e)$. 
Note that $\delta_e$ is a subset of $a$. 

Next, let us consider the case of $A(e)=\{(a,a),(b,b)\}$. 
We may assume that $a$ is the edge between $v_{i-1}$ and $v_i$, 
and that $b$ is the edge between $v_i$ and $v_{i+1}$. 
In this case we let $\delta_e$ be the edge $a$. 

By using these $\delta_e$'s, 
we would like to define continuous functions $\mu^e_n:X\to\{0,1\}$ 
for $e\in{\cal C}$ and $n\in\N$. 
Let $x=(x_n)_n\in X$. 
At first we define $\mu^e_1$ by 
\[ \mu^e_1(x)=\begin{cases}
1 & r(x_1)=r(e)\text{ and }
\lambda^{-1}(a+x_1)\subset\delta_e\text{ for some }a\in E(s(x_1)) \\
0 & \text{otherwise}. \end{cases} \]
For $n=2,3,\dots$, we define $\mu^e_n$ by 
\[ \mu^e_n(x)=\begin{cases}
\mu^e_{n-1}(x) & x_n=e \\
1 & x_n\neq e,r(x_n)=r(e)\text{ and }
\lambda^{-1}(a+x_n)\subset\delta_e\text{ for some }a\in E(s(x_n)) \\
0 & \text{otherwise}. \end{cases} \]
It is easily checked that $\mu^e_n$ is well-defined and continuous. 

The following is an immediate consequence of 
the definition of $\mu^e_n$. 

\begin{lem}\label{mu}
For $e\in{\cal C}$ and $(x,y)\in{\cal X}_n$, if 
\[ \mu^e_n(x)=\mu^e_n(y) \]
then we have 
\[ \mu^e_m(x)=\mu^e_m(y) \]
for all $m>n$. 
\end{lem}

For every $n\in\N$, 
we define a subset ${\cal X}'_n$ of ${\cal X}_n$ by 
\[ {\cal X}'_n=\{(x,y)\in{\cal X}_n:
\mu^e_n(x)=\mu^e_n(y)
\text{ for all }e\in{\cal C}\}. \]

\begin{lem}
For every $n\in\N$, 
${\cal X}'_n$ is a compact open subequivalence relation 
of ${\cal X}_n$, 
and ${\cal X}'_n$ is contained in ${\cal X}'_{n+1}$. 
\end{lem}
\begin{proof}
It is obvious that ${\cal X}'_n$ is a subequivalence relation 
of ${\cal X}_n$. 
Since $\mu^e_n$ is continuous, ${\cal X}'_n$ is compact and open. 
From the lemma above we can see ${\cal X}'_n\subset{\cal X}'_{n+1}$. 
\end{proof}

Define 
\[ {\cal X}'=\bigcup_{n\in\N}{\cal X}'_n. \]
By the lemma above, 
${\cal X}'$ is an AF equivalence relation on $X$. 

\begin{lem}\label{nosplit}
Suppose that $T\in\Omega_{punc}$ is not of type III. 
Then we have $[\pi(T)]_{\cal X}=[\pi(T)]_{{\cal X}'}$. 
\end{lem}
\begin{proof}
Let $x=(x_n)_n=\pi(T)$. 
Take $y=(y_n)_n\in[x]_{\cal X}$ arbitrarily. 
We would like to show that $y$ belongs to $[x]_{{\cal X}'}$. 
There exists $m\in\N$ such that $(x,y)\in{\cal X}_m$. 
Let $e\in{\cal C}$. 
If $x_l=e$ for every $l>m$, then 
$x$ is equivalent to $(e,e,\dots)$ in ${\cal X}$. 
It follows that $\pi^{-1}(x)=T$ is of type III. 
Hence we can find $l>m$ such that $x_l\neq e$. 
Combined with $x_l=y_l$, this implies $\mu^e_l(x)=\mu^e_l(y)$. 
By Lemma \ref{mu} the proof is completed. 
\end{proof}

\begin{lem}\label{minimal}
The equivalence relation ${\cal X}'$ is minimal. 
\end{lem}
\begin{proof}
Let $x=(x_n)_n\in X$. 
It suffices to show that $[x]_{{\cal X}'}$ is dense in $X$. 
If $\pi^{-1}(T)$ is not of type III, by the lemma above, 
$[x]_{{\cal X}'}=[x]_{\cal X}$. 
Since ${\cal X}$ is minimal, we can deduce that 
$[x]_{{\cal X}'}=[x]_{\cal X}$ is dense in $X$. 

Let us assume that $T=\pi^{-1}(x)$ is of type III. 
There exist $e\in{\cal C}$ and $m\in\N$ such that 
$x_n=e$ for all $n>m$. 
Take a finite path $(y_1,y_2,\dots,y_k)$ in $({\cal V},{\cal E})$ 
arbitrarily. 
It suffices to find $x'\in[x]_{{\cal X}'}$ 
whose initial segments agree with $(y_1,y_2,\dots,y_k)$. 
Since the directed graph $({\cal V},{\cal E})$ is primitive, 
we can find $N\in\N$ such that for any $p\in{\cal V}$ 
there exists a path of length $N$ which starts from $r(y_k)$ and 
ends at $p$. 
Let $l$ be a natural number greater than $m$ and $k+N$. 

We would like to consider the case of $\mu^e_l(x)=1$. 
We can find $f\in{\cal E}$ such that $f\neq e$, $r(f)=s(e)$ and 
$\lambda^{-1}(a+f)\subset\delta_e$ for some $a\in E(s(f))$. 
Take a finite path $(y_{k+1},y_{k+2},\dots,y_{l-1})$ 
in $({\cal V},{\cal E})$ so that $r(y_k)=s(y_{k+1})$ and 
$r(y_{l-1})=s(f)$. 
Define 
\[ x'=(y_1,y_2,\dots,y_{l-1},f,e,e,\dots)\in X. \]
By the choice of $f$ we have $\mu^e_l(x')=1$. 
Therefore we can conclude that $(x,x')\in {\cal X}'_{l+1}$. 

In the case of $\mu^e_l(x)=0$, 
we choose $f\in{\cal E}$ so that $r(f)=s(e)$ and 
$\lambda^{-1}(a+f)\subset\delta_e$ does not hold for any $a\in E(s(f))$. 
Then the proof goes in a similar fashion to the preceding paragraph. 
\end{proof}

It is well-known that both $R_{punc}$ and $R_{AF}$ 
are uniquely ergodic (see \cite{P} for example). 
It follows that both ${\cal X}_\omega$ and ${\cal X}$ are 
also uniquely ergodic. 
The next lemma claims that ${\cal X}'$ is uniquely ergodic, too. 

\begin{lem}\label{uniquely}
If $\nu$ is a ${\cal X}'$-invariant probability measure on $X$, then 
$\nu$ equals the unique ${\cal X}$-invariant probability measure. 
\end{lem}
\begin{proof}
It follows from Lemma \ref{minimal} that $\nu$ is nonatomic. 
By Lemma \ref{finiteIII} and Lemma \ref{nosplit}, 
\[ \{x\in X:[x]_{\cal X}\neq[x]_{{\cal X}'}\} \]
consists of finitely many ${\cal X}'$-orbits. 
Hence $\nu$ is ${\cal X}$-invariant. 
\end{proof}

\begin{lem}\label{empty}
Let $t_1,t_2\in\T$ be distinct elements. 
For 
\[ U=\{(x_n)_n\in X:\text{the cardinality of }A(x_1)\text{ is one}\}, \]
we have ${\cal X}'\cap((U\cap B_{t_1})\times(U\cap B_{t_2}))=\emptyset$. 
\end{lem}
\begin{proof}
For a proof by contradiction, assume that 
${\cal X}'\cap((U\cap B_{t_1})\times(U\cap B_{t_2}))$ contains $(x,y)$. 
Let $x=(x_n)_n$ and $y=(y_n)_n$. 
Since $t_1$ is distinct from $t_2$, $\pi^{-1}(x)$ is of type III. 
By Lemma \ref{periodicIII}, 
there exist $e\in{\cal C}$ and $m\in\N$ such that 
$x_n=y_n=e$ for all $n>m$. 
Besides, $A(e)$ should be $\{(a,a),(b,b)\}$ 
for some $a,b\in E(r(e))$. 
We may assume that $\delta_e=a$, $\theta_{r(e)}(a)=t_1$ and 
$\theta_{r(e)}(b)=t_2$. 
Then, for any $n>m$, one can see that 
$\mu^e_n(x)=1$ and $\mu^e_n(y)=0$ by an inductive calculation. 
But this contradicts $(x,y)\in{\cal X}'$. 
\end{proof}

\section{Affability}

Let $p\in{\cal V}$ be a prototile and 
let $a\in E(p)$ be an edge with $\theta_p(a)=t$. 
Thanks to the conditions (P5) and (P6), 
we can choose $e_i\in{\cal E}$ and $a_i\in E(s(e_i))$ for $i=1,2,3$ 
so that the following are satisfied. 
\begin{itemize}
\item $r(e_1)=s(e_2)$, $r(e_2)=s(e_3)$ and $r(e_3)=p$. 
\item $A(e_1)=\{(a_1,a_2)\}$, $A(e_2)=\{(a_2,a_3)\}$ and 
$A(e_3)=\{(a_3,a)\}$. 
\item $e_3$ does not belong to ${\cal C}$. 
\end{itemize}
Define a clopen subset $U_a$ of $X$ by 
\[ U_a=\{(x_n)_n\in X:x_1=e_1,x_2=e_2\text{ and }x_3=e_3\}. \]
Let us consider the prototile $s(e_3)$ and its edge $a_3$. 
Since the substitution system forces its border and 
satisfies the condition (P4), 
there exist $v\in\R^2$, $q\in{\cal V}$ and $b\in E(q)$ such that 
the following hold. 
\begin{itemize}
\item If $S\in\Omega$ contains $s(e_3)$, 
then $\omega(S)$ contains $q+v$. 
\item $\theta_q(b)=-t$ and $\lambda^{-1}(b+v)\subset a_3$. 
\end{itemize}
From the condition (P5), we can find $f\in{\cal E}$ such that 
$r(f)=q$ and $A(f)=\{(b',b)\}$ for some $b'\in s(f)$. 

Put $u_a=f+\lambda v-\lambda e_2-e_1\in\R^2$. 
For $T\in\pi^{-1}(U_a)$ we would like to consider $T-u_a$. 
By the definition of $U_a$, the tiling $\omega^{-2}(T)$ contains 
the tile $s(e_3)-\lambda^{-1}e_2-\lambda^{-2}e_1$, 
which implies 
\[ q+v-e_2-\lambda^{-1}e_1\in\omega^{-1}(T). \]
From $s(f)+f\in\omega(q)$ we get 
\[ s(f)+f+\lambda v-\lambda e_2-e_1=s(f)+u_a\in T. \]
Hence the map sending $\pi(T)\in U_a$ to $\pi(T-u_a)\in X$ is 
well-defined. 
We denote this map by $\beta_a:U_a\to X$. 

Let $T'=T-u_a$. 
The tile $s(f)\in T'$ contains the origin in its interior, 
The first coordinate of $\pi(T')$ is $f$ and $A(f)=\{(b',b)\}$. 
We define $w=-\lambda^{-1}e_3-\lambda^{-3}f-\lambda^{-2}v\in\R^2$. 
Then we have 
\begin{align*}
\omega^{-3}(T')
&= \omega^{-3}(T)-\lambda^{-3}u_a \\
&\ni r(e_3)-\lambda^{-1}e_3-\lambda^{-2}e_2-\lambda^{-3}e_1
-\lambda^{-3}u_a \\
&= p-\lambda^{-1}e_3-\lambda^{-3}f-\lambda^{-2}v \\
&= p+w. 
\end{align*}
Moreover we can see 
\begin{align*}
a+w
&\supset \lambda^{-1}(a_3+e_3)+w \\
&= \lambda^{-1}(a_3-\lambda^{-2}f-\lambda^{-1}v) \\
&\supset \lambda^{-1}(\lambda^{-1}(b+v)-\lambda^{-2}f-\lambda^{-1}v) \\
&=\lambda^{-2}(b-\lambda^{-1}f) \\
&\supset \lambda^{-2}(\lambda^{-1}(b'+f)-\lambda^{-1}f) \\
&= \lambda^{-3}b'. 
\end{align*}
We can summarize these arguments as follows. 

\begin{lem}\label{construction}
Let $p\in{\cal V}$ and $a\in E(p)$. Put $t=\theta_p(a)$. 
Define a clopen subset $U_a\subset X$ and 
a map $\beta_a:U_a\to X$ as above. 
\begin{enumerate}
\item $\beta_a$ is a homeomorphism from $U_a$ to $\beta_a(U_a)$ and 
$\{(x,\beta_a(x))\in X\times X:x\in U_a\}$ is 
a clopen subset of ${\cal X}_\omega$. 
\item If $x=(x_n)_n\in U_a$, then $x_3\notin{\cal C}$. 
\item If $x=(x_n)_n\in U_a$, then $A(x_1)=\{(a_1,a_2)\}$, 
$A(x_2)=\{(a_2,a_3)\}$ and $A(x_3)=\{(a_3,a)\}$
for some $a_1,a_2,a_3$ with $\theta_{s(x_1)}(a_1)=t$. 
\item If $y=(y_n)_n\in\beta_a(U_a)$, then $A(y_1)=\{(b',b)\}$ 
for some $b',b$ with $\theta_{s(y_1)}(b')=-t$. 
Moreover, there exists $w\in\R^2$ such that 
$p+w\in\omega^{-3}(\pi^{-1}(y))$ and $\lambda^{-3}b'\subset a+w$. 
\item For $x\in U_a$, $x\in B_t$ if and only if $\beta_a(x)\in B_{-t}$. 
\end{enumerate}
\end{lem}
\begin{proof}
Only (5) needs a proof. 
We use the notation used in the argument above. 
Take $x\in U_a$ and put $T=\pi^{-1}(x)$. 
We notice that the tiling $\omega^{-3}(T)$ contains the tile 
$p-\lambda^{-1}e_3-\lambda^{-2}e_2-\lambda^{-3}e_1$ and 
$a$ is an edge of $p$. 
Hence it is not hard to see that $x\in B_t$ if and only if 
\[ a-\lambda^{-1}e_3-\lambda^{-2}e_2-\lambda^{-3}e_1
\subset\lambda^n\partial(\omega^{-3-n}(T)) \]
for all $n\in\N$. 
We also notice that the tiling $T-u_a$ contains the tile $s(f)$ 
and $b'$ is an edge of $s(f)$. 
It follows that $\beta_a(x)\in B_{-t}$ if and only if 
\[ b'\subset\lambda^n\partial(\omega^{-n}(T-u_a)) \]
for all $n\in\N$. 
Therefore we have $\beta_a(x)\in B_{-t}$ if and only if 
\[ \lambda^{-3}(b'-u_a)\subset\lambda^n\partial(\omega^{-3-n}(T)) \]
for all $n\in\N$. 
From 
\[ \lambda^{-3}(b'-u_a)\subset 
a-\lambda^{-1}e_3-\lambda^{-2}e_2-\lambda^{-3}e_1 \]
and the condition (P1), 
we get $x\in B_t$ if and only if $\beta_a(x)\in B_{-t}$. 
\end{proof}

\begin{lem}\label{disjoint}
Let $p_1,p_2\in{\cal V}$ and let $a_i\in E(p_i)$ for $i=1,2$. 
\begin{enumerate}
\item If $a_1$ is not equal to $a_2$, then 
$U_{a_1}\cap U_{a_2}=\emptyset$. 
\item If $a_1$ is not equal to $a_2$, then 
$\beta_{a_1}(U_{a_1})\cap\beta_{a_2}(U_{a_2})=\emptyset$. 
\end{enumerate}
\end{lem}
\begin{proof}
(1) is evident from Lemma \ref{construction} (3). 
Let us prove (2). 
Suppose that $y=(y_n)_n\in X$ belongs to 
both $\beta_{a_1}(U_{a_1})$ and $\beta_{a_2}(U_{a_2})$. 
By Lemma \ref{construction} (4), 
we have $A(y_1)=\{(c,d)\}$ for some $c$ and $d$, 
and $-\theta_{s(y_1)}(c)$ must be equal to 
both $\theta_{p_1}(a_1)$ and $\theta_{p_2}(a_2)$. 
Thus, we can put 
$t=-\theta_{s(y_1)}(c)=\theta_{p_1}(a_1)=\theta_{p_2}(a_2)$. 
Besides, there exists $w_i\in\R^2$ such that 
$p_i+w_i\in\omega^{-3}(\pi^{-1}(y))$ and 
$\lambda^{-3}c\subset a_i+w_i$ for each $i=1,2$. 
Since $\omega^{-3}(\pi^{-1}(y))$ is a tiling, 
we can conclude $p_1=p_2$ and $a_1=a_2$. 
\end{proof}

Let $U_\omega$ be the union of all $U_a$'s, that is, 
\[ U_\omega=\bigcup_{a\in E(p),p\in{\cal V}}U_a. \]
For $x\in U_a\subset U_\omega$ we define $\bar{\beta}(x)=\beta_a(x)$. 
By the lemma above $\bar{\beta}$ is a well-defined homeomorphism 
from $U_\omega$ to $\bar{\beta}(U_\omega)$. 
Moreover $\{(x,\bar{\beta}(x)):x\in U_\omega\}$ is 
a clopen subset of ${\cal X}_\omega$. 
\bigskip

Choose a finite subset $F$ of $\T$ so that the following are satisfied: 
\begin{itemize}
\item If $B_t$ is not empty, then either of $t$ or $-t$ belongs to $F$. 
\item If $t\in F$, then $-t\notin F$. 
\end{itemize}
We define closed subsets $B$ and $B^*$ of $X$ by 
\[ B=U_\omega\cap\bigcup_{t\in F}B_t \]
and 
\[ B^*=\bar{\beta}(U_\omega)\cap\bigcup_{t\in F}B_{-t}. \]
(3) and (4) of Lemma \ref{construction} tells us that 
if $x\in U_a$, $a\in E(p)$ and $\theta_p(a)=t$, then 
$U_a\cap B_s=\emptyset$ and $\beta_a(U_a)\cap B_{-s}=\emptyset$ 
for all $s\neq t$. 
Moreover, from Lemma \ref{construction} (5), 
one can see $\bar{\beta}(B)=B^*$. 
We denote the restriction of $\bar{\beta}$ to $B$ by $\beta$. 

\begin{lem}\label{etalethin}
Both $B$ and $B^*$ are closed ${\cal X}'$-\'etale thin subsets. 
\end{lem}
\begin{proof}
At first let us prove that $B$ and $B^*$ are thin for ${\cal X}'$. 
Let $\nu$ be the unique ${\cal X}$-invariant probability measure 
on $X$. 
We remark that $\nu$ is also the unique ${\cal X}_\omega$-invariant 
probability measure. 
By Lemma \ref{uniquely} it suffices to show $\nu(B\cup B^*)=0$. 
But this follows immediately from \cite[Theorem 2.1]{P}. 

We next show that $B$ is \'etale for ${\cal X}'$. 
Note that $B$ is equal to 
a disjoint union $\bigcup_{t\in F}U_\omega\cap B_t$. 
Suppose that $(x,y)$ belongs to ${\cal X}'\cap(B\times B)$. 
Let $x\in U_\omega\cap B_{t_1}$ and $y\in U_\omega\cap B_{t_2}$. 
By Lemma \ref{empty} and Lemma \ref{construction} (3), 
$t_1$ must equal $t_2$. 
Thus, $x$ and $y$ lie in some $U_\omega\cap B_t$. 
From Proposition \ref{etale}, 
$U_\omega\cap B_t$ is ${\cal X}$-\'etale. 
Since ${\cal X}'$ is an open subequivalence relation of ${\cal X}$, 
we can deduce that $B$ is ${\cal X}'$-\'etale. 

Similarly we can prove that $B^*$ is also ${\cal X}'$-\'etale. 
The proof is completed. 
\end{proof}

\begin{lem}\label{norelation}
We have ${\cal X}'\cap(B\times B^*)=\emptyset$. 
\end{lem}
\begin{proof}
For a proof by contradiction, suppose that 
${\cal X}'\cap(B\times B^*)$ contains $(x,y)$. 
Let $x\in U_\omega\cap B_{t_1}$ and 
$y\in \beta(U_\omega)\cap B_{-t_2}$, 
where $t_1$ and $t_2$ belong to $F$. 
By Lemma \ref{construction} (3), (4) and Lemma \ref{empty}, 
$t_1$ must equal $-t_2$. 
But this contradicts the choice of $F$. 
\end{proof}

\begin{lem}\label{iso}
The homeomorphism $\beta:B\to B^*$ induces an isomorphism 
between ${\cal X}'\cap(B\cap B)$ and ${\cal X}'\cap(B^*\times B^*)$. 
\end{lem}
\begin{proof}
Since $\{(x,\bar{\beta}(x)):x\in U_\omega\}$ is 
a clopen subset of ${\cal X}_\omega$, 
\[ \bar{\beta}\times\bar{\beta}:
{\cal X}_\omega\cap(U_\omega\cap U_\omega)\to
{\cal X}_\omega\cap(\bar{\beta}(U_\omega)\times\bar{\beta}(U_\omega)) \]
is a well-defined isomorphism. 
Because ${\cal X}'$ is an open subequivalence relation of 
${\cal X}_\omega$ and 
$\beta$ is a restriction of 
$\bar{\beta}:U_\omega\to\bar{\beta}(U_\omega)$ to $B$, 
it suffices to show that 
$\beta$ induces a bijective correspondence 
between ${\cal X}'\cap(B\times B)$ and ${\cal X}'\cap(B^*\times B^*)$. 

Take $x,y\in B$. 
It suffices to show that $(x,y)\in{\cal X}'$ 
if and only if $(\beta(x),\beta(y))\in{\cal X}'$. 
When $\pi^{-1}(x)$ is of type I, 
we have nothing to do because of $[x]_{{\cal X}_\omega}=[x]_{{\cal X}'}$. 

Assume that $T=\pi^{-1}(x)$ is of type II. 
In this case $\partial_\infty(T)$ is a line and 
it divides the plane into two open half-planes. 
By Lemma \ref{nosplit} we have $[x]_{\cal X}=[x]_{{\cal X}'}$, 
and so it suffices to show $(x,y)\in{\cal X}$ 
if and only if $(\beta(x),\beta(y))\in{\cal X}$. 
If $(x,y)\in{\cal X}$, then there exists $t\in F$ such that 
$x,y\in B_t$. 
It follows that $\beta(x),\beta(y)\in B_{-t}$, 
which means $(\beta(x),\beta(y))\in{\cal X}$. 
In the same way we can show that 
$(\beta(x),\beta(y))\in{\cal X}$ implies $(x,y)\in{\cal X}$. 

Suppose that $T=\pi^{-1}(x)$ is of type III. 
Thus $\partial_\infty(T)$ consists of finitely many half lines. 
Let $v\in\R^2$ be the center of $\partial_\infty(T)$. 
Note that if $p\in T$ does not meet $v$, then 
$p$ meets at most one half line in $\partial_\infty(T)$. 
We denote the canonical bijective correspondence 
from $[x]_{{\cal X}_\omega}$ to $T$ by $\rho$, that is, 
if the puncture of $p\in T$ is $u\in\R^2$ then $\rho(\pi(T-u))=p$. 
Suppose that $(x,y)$ belongs to ${\cal X}'$. 
Let $p=\rho(x)$ and $q=\rho(y)$. 
By Lemma \ref{construction} (2), 
we notice that neither $p$ nor $q$ meets $v$. 
From $(x,y)\in{\cal X}'$ we can deduce that 
there exists a half line $\ell$ in $\partial_\infty(T)$ such that 
both $p$ and $q$ meet $\ell$ from the same side of $\ell$. 
It follows that both $\rho(\beta(x))$ and $\rho(\beta(y))$ meet $\ell$ 
but they lie on the opposite side from $p$ and $q$ 
across the half line $\ell$. 
Hence we can see that $(\beta(x),\beta(y))\in{\cal X}'$. 
Next, assume that $(x,y)$ belongs to ${\cal X}_\omega$ 
but does not belong to ${\cal X}'$. 
There exist half lines $\ell_1$ and $\ell_2$ in $\partial_\infty(T)$ 
such that $p$ meets $\ell_1$ and $q$ meets $\ell_2$. 
If $\ell_1$=$\ell_2$, then 
$p$ and $q$ never lie on the same side of $\ell_1=\ell_2$. 
But, if $p$ is located on the opposite side from $q$ 
across the half line $\ell_1=\ell_2$, then 
$F\subset\T$ must contain both $t$ and $-t$ for some $t\in\T$, 
which is a contradiction. 
Therefore $\ell_1$ and $\ell_2$ should be different. 
Since $\rho(\beta(x))$ meets $\ell_1$ and $\rho(\beta(y))$ meets $\ell_2$, 
we can conclude that $(\beta(x),\beta(y))$ never belongs to ${\cal X}'$. 
Consequently we have $(x,y)\in{\cal X}'$ if and only if 
$(\beta(x),\beta(y))\in{\cal X}'$, 
thereby completing the proof. 
\end{proof}

Let $\widetilde{\cal X}$ be the equivalence relation 
generated by ${\cal X}'$ and $\{(x,\beta(x))\in X\times X:x\in B\}$. 
It is obvious that $\widetilde{\cal X}$ is 
a subequivalence relation of ${\cal X}_\omega$. 

\begin{lem}\label{absorp}
The equivalence relation $\widetilde{\cal X}$ is 
orbit equivalent to ${\cal X}$. 
\end{lem}
\begin{proof}
From Lemma \ref{minimal}, \ref{etalethin}, \ref{norelation}, \ref{iso}, 
we can apply the absorption theorem \cite[Theorem 4.18]{GPS2} 
to ${\cal X}'$ and $\beta:B\to B^*$. 
It follows that $\widetilde{\cal X}$ is orbit equivalent to ${\cal X}'$. 
Both ${\cal X}'$ and ${\cal X}$ are minimal AF equivalence relations. 
Moreover they are uniquely ergodic 
with the same invariant probability measure $\nu$ 
by Lemma \ref{uniquely}. 
By \cite[Theorem 2.3]{GPS1}, we can conclude that 
$\widetilde{\cal X}$ is orbit equivalent to ${\cal X}$. 
\end{proof}

Now we are ready to prove the main theorem. 

\begin{thm}
Let $({\cal V},\omega)$ be a substitution tiling system in $\R^2$ 
which is primitive, aperiodic and satisfies the finite pattern condition. 
Suppose that each prototile is a polygon and 
the conditions (P1), (P2) and (P3) are satisfied. 
Then the equivalence relation $R_{punc}$ on $\Omega_{punc}$ is 
orbit equivalent to $R_{AF}$. 
In particular $R_{punc}$ is affable. 
\end{thm}
\begin{proof}
Let $x=\pi(T)\in X$. 
When $T$ is of type I, 
the ${\cal X}_\omega$-equivalence class $[x]_{{\cal X}_\omega}$ is equal to 
the ${\cal X}'$-equivalence class $[x]_{{\cal X}'}$. 
Hence $[x]_{{\cal X}_\omega}=[x]_{\widetilde{\cal X}}$. 
Let us assume that $T$ is of type II. 
The ${\cal X}_\omega$-orbit $[x]_{{\cal X}_\omega}$ splits 
into two ${\cal X}'$-orbits. 
As described in the proof of Lemma \ref{iso}, 
these two orbits are glued by $\beta:B\to B^*$. 
It follows that $[x]_{{\cal X}_\omega}$ agrees 
with $[x]_{\widetilde{\cal X}}$. 
If $T$ is of type III, then 
$[x]_{{\cal X}_\omega}$ may not agree with $[x]_{\widetilde{\cal X}}$. 
By Lemma \ref{finiteIII}, however, 
there are only finitely many $R_{punc}$-orbits of type III. 
In addition, every ${\cal X}_\omega$-orbit splits into 
finitely many ${\cal X}'$-orbits. 
Therefore, we can find $(x_1,y_1),(x_2,y_2),\dots,(x_n,y_n)$ in $X\times X$ 
such that ${\cal X}_\omega$ is generated by $\widetilde{\cal X}$ and 
$(x_1,y_1),(x_2,y_2),\dots,(x_n,y_n)$. 
It follows from Lemma \ref{absorp} and \cite[Corollary 4.17]{GPS2} that 
${\cal X}_\omega$ is orbit equivalent to ${\cal X}$. 
Consequently the equivalence relation $R_{punc}$ on $\Omega_{punc}$ is 
orbit equivalent to $R_{AF}$. 
\end{proof}

\flushleft{
\textit{e-mail: matui@math.s.chiba-u.ac.jp \\
Graduate School of Science and Technology,\\
Chiba University,\\
1-33 Yayoi-cho, Inage-ku,\\
Chiba 263-8522,\\
Japan. }}

\end{document}